\documentclass[11pt,reqno]{amsart}
\usepackage{amscd,amssymb,amsmath,amsthm}
\usepackage[arrow,matrix]{xy}
\usepackage{graphicx}
\usepackage{cite}
\newtheorem{thm}[subsection]{Theorem}
\newtheorem{lemma}[subsection]{Lemma}
\newtheorem{pro}[subsection]{Proposition}
\newtheorem{cor}[subsection]{Corollary}

\newtheorem{defn}[subsection]{Definition}

\numberwithin{equation}{section} 

\newcommand{\s}{{\sigma}}
\newcommand{\de}{{\xi}}

\newcommand{\bea}{\begin{eqnarray}}
\newcommand{\eea}{\end{eqnarray}}

\begin{document}
\title[Non-unique Gibbs measure of a model]{Non-uniqueness of Gibbs Measure for Models
With Uncountable Set of Spin Values on a Cayley Tree}

\author{Yu. Kh. Eshkabilov, F. H. Haydarov, U. A. Rozikov}

 \address{Yu.\ Kh.\ Eshkabilov\\ National University of Uzbekistan,
Tashkent, Uzbekistan.}
\email {yusup62@mail.ru}

\address{F.\ H.\ Haydarov\\ National University of Uzbekistan,
Tashkent, Uzbekistan.}
\email{haydarov$_-$imc@mail.ru}

 \address{U.\ A.\ Rozikov\\ Institute of mathematics and information technologies,
Tashkent, Uzbekistan.}
\email {rozikovu@yandex.ru}

\begin{abstract} In this paper we construct several models with nearest-neighbor
interactions and with the set $[0,1]$ of spin values, on a Cayley
tree of order $k\geq 2$. We prove that each of the constructed model
has at least two translational-invariant Gibbs measures.
\end{abstract}
\maketitle

{\bf Mathematics Subject Classifications (2010).} 82B05, 82B20 (primary);
60K35 (secondary)

{\bf{Key words.}} Cayley tree, configuration, Gibbs measures, uniqueness.

\section{Introduction} \label{sec:intro}

Spin systems on lattices are a large class of systems considered in statistical
mechanics. Some of them have a real physical meaning, others are studied as
suitably simplified models of more complicated systems.
The structure of the lattice (graph) plays an important role in investigations
of spin systems. For example, in order to study the phase transition problem
for a system on $Z^d$ and on Cayley
tree there are two different methods: Pirogov-Sinai theory on $Z^d$, Markov
random field theory and recurrent equations of this theory on Cayley tree.
In \cite{2}-\cite{4},\cite{p1a}, \cite{11}-\cite{p2},
\cite{13}, \cite{14}, \cite{16} for several models on Cayley tree, using the Markov random field theory
Gibbs measures are described.

These papers are devoted to models with a {\it finite} set of spin values.
Mainly were shown that these models have finitely many translation-invariant  and uncountable
numbers  of the  non-translation-invariant extreme Gibbs
measures. Also for several models (see, for example, \cite{p1, p1a, p2}) it were proved that there exist
three  periodic  Gibbs  measures (which are invariant with respect to normal  subgroups  of  finite index of the group representation of the Cayley tree) and there are
uncountable number of non-periodic Gibbs measures.

In \cite{6} the Potts model with
a {\it countable} set of spin values on a Cayley tree is considered and it was showed
that the set of translation-invariant splitting Gibbs measures of the model contains at most one
point, independently on parameters of the Potts model with countable set of spin values
on the Cayley tree. This is a crucial difference from the models with a finite set of spin values, since the last ones may have more than one translation-invariant Gibbs measures.

This paper is continuation of our investigations \cite{ehr},\cite{re}.
In \cite{re} models (Hamiltonians) with nearest-neighbor interactions and with the
({\it uncountable}) set $[0,1]$ of spin values, on a Cayley tree of order $k\geq 1$ were studied.

A central problem in the theory of Gibbs measures is to
describe infinite-volume (or limiting) Gibbs measures
corresponding to a given Hamiltonian. In \cite{re} we reduced the problem  to the description
of the solutions of some nonlinear integral equation. Then for $k = 1$ we showed that
the integral equation has a unique solution. In case $k\geq 2$ some models (with
the set $[0, 1]$ of spin values) which have a unique splitting Gibbs measure are
constructed.  In our next paper \cite{ehr} it was found a sufficient condition on
Hamiltonian of the model with an uncountable set of spin values under which the model has unique
translation-invariant splitting Gibbs measure. But we had not any example of model
(with uncountable spin values) with more than one translation-invariant Gibbs measure.

This problem is solved in this paper: we shall construct several models with nearest-neighbor
interactions and with the set $[0,1]$ of spin values, on a Cayley
tree of order $k\geq 2$. We prove that each of the constructed model
have at least two translational-invariant Gibbs measures.

\section{Preliminaries}

A Cayley tree $\Gamma^k=(V,L)$ of order $k\geq 1$ is an infinite
homogeneous tree, i.e., a graph without cycles, with
exactly $k+1$ edges incident to each vertices. Here $V$ is the set
of vertices and $L$ that of edges (arcs).

Consider models where the spin takes values in the set $[0,1]$, and is assigned to the vertexes
of the tree. For $A\subset V$ a configuration $\s_A$ on $A$ is an arbitrary function $\s_A:A\to
[0,1]$. Denote $\Omega_A=[0,1]^A$ the set of all configurations on $A$. A configuration $\sigma$ on
$V$ is then defined as a function $x\in V\mapsto\sigma (x)\in [0,1]$; the set of all configurations
is $[0,1]^V$. The (formal) Hamiltonian of the model is :
\begin{equation}\label{e1.1}
 H(\sigma)=-J\sum_{\langle x,y\rangle\in L}
\de_{\sigma(x)\sigma(y)},
\end{equation}
where $J \in R\setminus \{0\}$
and $\de: (u,v)\in [0,1]^2\to \de_{uv}\in R$ is a given bounded,
measurable function. As usually, $\langle x,y\rangle$ stands for
nearest neighbor vertices.

Let $\lambda$ be the Lebesgue measure on $[0,1]$.  On the set of all
configurations on $A$ the a priori measure $\lambda_A$ is introduced as
the $|A|$fold product of the measure $\lambda$. Here and further on
$|A|$ denotes the cardinality of $A$.   We consider a standard
sigma-algebra ${\mathcal B}$ of subsets of $\Omega=[0,1]^V$ generated
by the measurable cylinder subsets.
 A probability measure $\mu$ on $(\Omega,{\mathcal B})$
is called a Gibbs measure (with Hamiltonian $H$) if it satisfies the DLR equation, namely for any
$n=1,2,\ldots$ and $\sigma_n\in\Omega_{V_n}$:
$$\mu\left(\left\{\sigma\in\Omega :\;
\sigma\big|_{V_n}=\sigma_n\right\}\right)= \int_{\Omega}\mu ({\rm
d}\omega)\nu^{V_n}_{\omega|_{W_{n+1}}} (\sigma_n),$$ where $\nu^{V_n}_{\omega|_{W_{n+1}}}$ is the
conditional Gibbs density
$$ \nu^{V_n}_{\omega|_{W_{n+1}}}(\sigma_n)=\frac{1}{Z_n\left(
\omega\big|_{W_{n+1}}\right)}\exp\;\left(-\beta H
\left(\sigma_n\,||\,\omega\big|_{W_{n+1}}\right)\right),
$$
and $\beta={1\over T}$, $T>0 $ is temperature.
Here and below, $W_l$ stands for a `sphere' and $V_l$ for a
`ball' on the tree, of radius $l=1,2,\ldots$,
centered at a fixed vertex $x^0$ (an origin):
$$W_l=\{x\in V: d(x,x^0)=l\},\;\;V_l=\{x\in V: d(x,x^0)\leq l\};$$
and
$$L_n=\{\langle x,y\rangle\in L: x,y\in V_n\};$$
distance $d(x,y)$, $x,y\in V$, is the length of (i.e. the number of edges in) the shortest path
connecting $x$ with $y$. $\Omega_{V_n}$ is the set of configurations in $V_n$ (and $\Omega_{W_n}$
that in $W_n$; see below). Furthermore, $\sigma\big|_{V_n}$ and $\omega\big|_{W_{n+1}}$ denote the
restrictions of configurations $\sigma,\omega\in\Omega$ to $V_n$ and $W_{n+1}$, respectively. Next,
$\sigma_n:\;x\in V_n\mapsto \sigma_n(x)$ is a configuration in $V_n$ and
$H\left(\sigma_n\,||\,\omega\big|_{W_{n+1}}\right)$ is defined as the sum
$H\left(\sigma_n\right)+U\left(\sigma_n, \omega\big|_{W_{n+1}}\right)$ where
$$H\left(\sigma_n\right)
=-J\sum_{\langle x,y\rangle\in L_n}\de_{\sigma_n(x)\sigma_n(y)},$$
$$U\left(\sigma_n,
\omega\big|_{W_{n+1}}\right)=
-J\sum_{\langle x,y\rangle:\;x\in V_n, y\in W_{n+1}}
\de_{\sigma_n(x)\omega (y)}.$$
Finally, $Z_n\left(\omega\big|_{W_{n+1}}\right)$
stands for the partition function in $V_n$, with
the boundary condition $\omega\big|_{W_{n+1}}$:
$$Z_n\left(\omega\big|_{W_{n+1}}\right)=
\int_{\Omega_{V_n}} \exp\;\left(-\beta H \left({\widetilde\sigma}_n\,||\,\omega
\big|_{W_{n+1}}\right)\right)\lambda_{V_n}(d{\widetilde\sigma}_n).$$

Due to the nearest-neighbor character of the interaction, the
Gibbs measure possesses a natural Markov property: for given a
configuration $\omega_n$ on $W_n$, random configurations in
$V_{n-1}$ (i.e., `inside' $W_n$) and in $V\setminus V_{n+1}$
(i.e., `outside' $W_n$) are conditionally independent.

We use a standard definition of a translation-invariant measure
(see, e.g., \cite{12}).
 The main object of study in this
paper are translation-invariant Gibbs measures for the model (\ref{e1.1})
on Cayley tree. In \cite{re} this problem of description of such
measures was reduced to the description of the solutions of a nonlinear
integral equation. For finite and countable sets of spin values
this argument is well known (see, e.g.
\cite{2}-\cite{6},\cite{11},\cite{13},\cite{14},\cite{16}).

Write $x<y$ if the path from $x^0$ to $y$ goes through $x$. Call vertex $y$ a direct successor of
$x$ if $y>x$ and $x,y$ are nearest neighbors. Denote by $S(x)$ the set of direct successors of $x$.
Observe that any vertex $x\ne x^0$ has $k$ direct successors and $x^0$ has $k+1$.

Let $h:\;x\in V\mapsto h_x=(h_{t,x}, t\in [0,1]) \in R^{[0,1]}$ be mapping of $x\in V\setminus
\{x^0\}$.  Given
$n=1,2,\ldots$, consider the probability distribution $\mu^{(n)}$ on $\Omega_{V_n}$ defined by
\begin{equation}\label{e2}
\mu^{(n)}(\sigma_n)=Z_n^{-1}\exp\left(-\beta H(\sigma_n)
+\sum_{x\in W_n}h_{\sigma(x),x}\right),
\end{equation}
 Here, as before, $\sigma_n:x\in V_n\mapsto
\sigma(x)$ and $Z_n$ is the corresponding partition function:
\begin{equation}\label{e3} Z_n=\int_{\Omega_{V_n}}
\exp\left(-\beta H({\widetilde\sigma}_n) +\sum_{x\in W_n}h_{{\widetilde\sigma}(x),x}\right)
\lambda_{V_n}({d\widetilde\s_n}).
\end{equation}

The probability distributions $\mu^{(n)}$ are compatible if for any $n\geq 1$ and
$\sigma_{n-1}\in\Omega_{V_{n-1}}$:
\begin{equation}\label{e4}
\int_{\Omega_{W_n}}\mu^{(n)}(\sigma_{n-1}\vee\omega_n)\lambda_{W_n}(d(\omega_n))=
\mu^{(n-1)}(\sigma_{n-1}).
\end{equation} Here
$\sigma_{n-1}\vee\omega_n\in\Omega_{V_n}$ is the concatenation of
$\sigma_{n-1}$ and $\omega_n$. In this case there exists a unique
measure $\mu$ on $\Omega_V$ such that, for any $n$ and
$\sigma_n\in\Omega_{V_n}$, $\mu \left(\left\{\sigma
\Big|_{V_n}=\sigma_n\right\}\right)=\mu^{(n)}(\sigma_n)$.

\begin{defn} The measure $\mu$ is called {\it splitting
Gibbs measure} corresponding to Hamiltonian (\ref{e1.1}) and function
$x\mapsto h_x$, $x\neq x^0$. \end{defn}

 The following
statement describes conditions on $h_x$ guaranteeing compatibility
of the corresponding distributions $\mu^{(n)}(\sigma_n).$

 \begin{pro}\label{p1}\cite{re} {\it The probability distributions
$\mu^{(n)}(\sigma_n)$, $n=1,2,\ldots$, in}
(\ref{e2}) {\sl are compatible iff for any $x\in V\setminus\{x^0\}$
the following equation holds:
\begin{equation}\label{e5}
 f(t,x)=\prod_{y\in S(x)}{\int_0^1\exp(J\beta\de_{tu})f(u,y)du \over \int_0^1\exp(J\beta{\de_{0u}})f(u,y)du}.
 \end{equation}
Here, and below  $f(t,x)=\exp(h_{t,x}-h_{0,x}), \ t\in [0,1]$ and
$du=\lambda(du)$ is the Lebesgue measure.}
\end{pro}

From Proposition \ref{p1} it follows that for any $h=\{h_x\in R^{[0,1]},\
\ x\in V\}$ satisfying (\ref{e5}) there exists a unique Gibbs measure
$\mu$ and vice versa. However, the analysis of solutions to (\ref{e5}) is
not easy. This difficulty depends on the given function $\xi$.

Let $\xi_{tu}$ is a continuous function
and we are going to construct functions $\xi_{tu}$ under which the equation (\ref{e5})
has at least two solutions in the class of
translational-invariant functions $f(t,x)$, i.e $f(t,x)=f(t),$ for
any $x\in V$. For such functions equation (\ref{e5}) can be written as
\begin{equation}\label{e1.2}
f(t)=\left({\int_0^1K(t,u)f(u)du\over \int_0^1 K(0,u)f(u)du}\right)^k,
\end{equation}
where $K(t,u)=\exp(J\beta \xi_{tu}), f(t)>0, t,u\in [0,1].$

We put
$$C^+[0,1]=\{f\in C[0,1]: f(x)\geq 0\}.$$
We are interested to positive continuous solutions to (\ref{e1.2}), i.e. such that

$f\in C_0^+[0,1]=\{f\in C[0,1]: f(x)\geq 0\}\setminus \{\theta\equiv 0\}$.

Note that equation (\ref{e1.2}) is not linear for any $k\geq 1$.

 Define the
operator $R_{k}:C^{+}_{0}[0,1]\rightarrow C^{+}_{0}[0,1]$ by
$$(R_{k}f)(t)=\left[\frac{(Wf)(t)}{(Wf)(0)}\right]^{k}, \,\
k\in\mathbb{N},$$\\
where $W:C[0,1]\rightarrow C[0,1]$ is linear operator, which is
defined by :
 $$(Wf)(t)=\int^{1}_{0}K(t,u)f(u)du.$$

Then the equation (\ref{e1.2}) can be written as

$$R_{k}f=f, \,\ f\in C^{+}_{0}[0,1].$$

\section{The Hammerstein's nonlinear integral equation}

For every $k\in\mathbb{N}$ we consider an integral operator
$H_{k}$ acting in the cone $C^{+}[0,1]$ as
$$(H_{k}f)(t)=\int^{1}_{0}K(t,u)f^{k}(u)du, \,\ k\in\mathbb{N}.$$

The operator $H_{k}$ is called Hammerstein's integral operator of
order $k$. Clearly that, if $k\geq2$ then $H_{k}$ is a
nonlinear operator.

\begin{lemma}\label{l2.1.} Let $k\geq2$. The equation
\begin{equation}\label{e2.1}
R_{k}f=f, \,\ f\in C^{+}_{0}[0,1]
\end{equation}
has a nontrivial positive solution iff the Hammerstein's operator
has a positive eigenvalue, i.e. the Hammerstein's equation
\begin{equation}\label{e2.2} H_{k}f=\lambda f, \,\ f\in C^{+}[0,1]
\end{equation}
has a nonzero positive solution for some $\lambda>0$.
\end{lemma}
\proof {\it Necessariness}. Let $f_{0}\in
C^{+}_{0}[0,1]$ be a solution of the equation (\ref{e2.1}). We have
$$(Wf_{0})(t)=(Wf_{0})(0)f_{0}^{\frac{1}{k}}(t).$$
From this equality we get
$$(H_{k}h)(t)=\lambda_{0}h(t),$$
where
$h(t)=\sqrt[k]{f_{0}(t)}\in C^{+}_{0}[0,1]$ and
$\lambda_{0}=(Wf_{0})(0),$ i.e., the number $\lambda_{0}$ is the
positive eigenvalue of the Hammerstein's operator $H_{k}.$

{\it Sufficiency}. Let $\lambda_{0}$ be a positive eigenvalue of the
operator $H_{k}.$ Then $\lambda_{0}>0$ and there exists $f_{0}\in
C^{+}_{0}[0,1]$ such that $$H_{k}f_{0}=\lambda_{0}f_{0}.$$
Obviously, the function $f_{0}(t)$ is a strictly positive. Put
$$f(t)=\frac{f_{0}(t)}{f_{0}(0)}, \,\ t\in[0,1].$$
Then the number $\lambda=\lambda_{0}f^{1-k}_{0}(0)$ is an eigenvalue of $H_{k}$ and
 corresponding the positive eigenfunction $f(t)$ satisfies

$$H_{k}f(t)=\frac{1}{f^{k}_{0}(0)}(H_{k}f_{0})(t)=\frac{\lambda_{0}}{f^{k}_{0}(0)}f_{0}(t)=\lambda_{0}f^{1-k}_{0}(0)f(t)=\lambda f(t). $$
Define
$$h(t)=\left(\frac{f_{0}(t)}{f_{0}(0)}\right)^{k}=f^{k}(t).$$
Then
$$(R_{k}h)(t)=\left(\frac{(Wh)(t)}{(Wh)(0)}\right)^{k}=\left(\frac{(H_{k}f)(t)}{(H_{k}f)(0)}\right)^{k}=
\left(\frac{\lambda_{0}f^{1-k}_{0}(0)f(t)}{\lambda_{0}f^{1-k}_{0}(0)f(0)}\right)^{k}=f^{k}(t)=h(t).$$
\endproof
\begin{cor}\label{c2.2.} Let $k\geq 2$. If a function $f\in
C^{+}_{0}[0,1]$ is an eigenfunction of $H_{k},$ then the function
$$h(t)=\left(\frac{f(t)}{f(0)}\right)^{k}$$
is a solution to the equation (\ref{e2.1}).
\end{cor}

\section{Existence of two Gibbs measures for the model (\ref{e1.1}): case $k=2$}

Consider the case $k=2$ in the model (\ref{e1.1}) and
$$\xi_{t,u}=\frac{1}{\beta J}\ln\left(1+\frac{14}{15}\cdot\sqrt[5]{4\left(t-\frac{1}{2}\right)\left(u-\frac{1}{2}\right)}\right), \,\ t,u\in[0,1].$$
Then, for the kernel $K(t,u)$ of the Hammerstein's integral
operator $H_{2}$ we have
$$K(t,u)=1+\frac{14}{15}\cdot\sqrt[5]{4\left(t-\frac{1}{2}\right)\left(u-\frac{1}{2}\right)}.$$

\begin{pro}\label{p3.1.} The Hammerstein's operator
$H_{2}:$
$$(H_{2}f)(t)=\int^{1}_{0}K(t,u)f^{2}(u)du$$
in the space $C[0,1]$ has at least two strictly positive fixed
points.
\end{pro}
\proof a) Let $f_{1}(t)\equiv1.$ Then we have
$$(H_{2}f_{1})(t)=1+\frac{14}{15}\cdot\sqrt[5]{4\left(t-\frac{1}{2}\right)}\cdot\int^{1}_{0}\left(u-
\frac{1}{2}\right)^{\frac{1}{5}}du=
1=f_{1}(t), \,\ t\in[0,1].$$
b) Denote
$$f_{2}(t)=\frac{3}{4}+\sqrt{\frac{21}{5}}\cdot\frac{\sqrt[5]{2}}{4}\cdot\left(t-\frac{1}{2}\right)^{\frac{1}{5}}, \,\ t\in[0,1].$$

Then $f_{2}\in C[0,1]$ and the function $f_{2}(t)$ is strictly
positive. Put
$$a=\frac{14}{15}\cdot\sqrt[5]{4}  , \,\ b=\sqrt{\frac{21}{5}}\cdot\frac{\sqrt[5]{2}}{4}.$$
We have $$H_{2}f_{2}=h_{1}(t)+h_{2}(t)+h_{3}(t)+\gamma ,$$
where

$$h_{1}(t)=ab^{2}\cdot\sqrt[5]{t-\frac{1}{2}}\cdot\int^{1}_{0}\sqrt[5]{\left(u-\frac{1}{2}\right)^{3}}du,
$$

$$h_{2}(t)=\frac{3ab}{2}\cdot\sqrt[5]{t-\frac{1}{2}}\cdot\int^{1}_{0}\sqrt[5]{\left(u-\frac{1}{2}\right)^{2}}du,$$

$$h_{3}(t)=\frac{9a}{16}\cdot\sqrt[5]{t-\frac{1}{2}}\cdot\int^{1}_{0}\sqrt[5]{u-\frac{1}{2}}
du,$$

$$\gamma=\int^{1}_{0}f^{2}_{2}(u)du.$$

It is clear that $$h_{1}(t)=h_{3}(t)\equiv0.$$ For the function
$h_{2}(t)$ we obtain
$$h_{2}(t)=\frac{3ab}{2}\cdot\sqrt[5]{t-\frac{1}{2}}\cdot\int^{1/2}_{-1/2}u^{\frac{2}{5}}du=\frac{15ab}{14\sqrt[5]{4}}\cdot\sqrt[5]{t-\frac{1}{2}} \,\ .$$
Observe that $$\gamma=\frac{5b^{2}}{7\sqrt[5]{4}}+\frac{9}{16}.$$
Consequently we have
$$(H_{2}f_{2})(t)=h_{2}(t)+\gamma=\frac{15ab}{14\sqrt[5]{4}}\cdot\sqrt[5]{t-\frac{1}{2}}+\frac{5b^{2}}{7\sqrt[5]{4}}+\frac{9}{16}=
\sqrt{\frac{21}{5}}\cdot\frac{\sqrt[5]{2}}{4}\cdot\sqrt[5]{t-\frac{1}{2}}+\frac{3}{4}=f_{2}(t).$$
\endproof
Denote by $\mu_1$ and $\mu_2$ the translation-invariant Gibbs measures which by Proposition \ref{p1} correspond to
solutions $f_1(t)=1$ and $f_{2}(t)=\frac{3}{4}+\sqrt{\frac{21}{5}}\cdot\frac{\sqrt[5]{2}}{4}\cdot\left(t-\frac{1}{2}\right)^{\frac{1}{5}}$.

Thus we have proved the following

\begin{thm}\label{t3.2.}  The model
$$H(\sigma)=-\frac{1}{\beta}\sum\limits_{<x,y>\atop{x,y}\in V}\ln\left(1+\frac{14}{15}\cdot
\sqrt[5]{4\left(\sigma(x)-\frac{1}{2}\right)\left(\sigma(y)-\frac{1}{2}\right)}
\right), \,\ \sigma\in\Omega_{V}$$ on the Cayley tree $\Gamma^2$
has at least two translation-invariant Gibbs measures $\mu_1$, $\mu_2$.
\end{thm}

\section{Existence of two Gibbs measures for the model (\ref{e1.1}): case $k=3$}

Now we shall consider the case $k=3$ and
$$\xi_{t,u}=\frac{1}{\beta J}\ln\left(1+\frac{1}{2}\cdot\sqrt[7]{4\left(t-\frac{1}{2}\right)\left(u-\frac{1}{2}\right)}\right), \,\ t,u\in[0,1].$$
Then, for the kernel $K(t,u)$ of the operator $H_{3}$ we have
$$K(t,u)=1+\frac{1}{2}\sqrt[7]{4\left(t-\frac{1}{2}\right)\left(u-\frac{1}{2}\right)}.$$

\begin{pro}\label{p4.1.} The operator $H_{3}:$
$$(H_3f)(t)=\int^{1}_{0}\left(1+\frac{1}{2}\cdot\sqrt[7]{4\left(t-\frac{1}{2}\right)\left(u-\frac{1}{2}\right)}\right)f^{3}(u)du$$
in the space $C[0,1]$ has at least two strictly positive fixed
points.
\end{pro}
\proof a) Let $f_{1}(t)\equiv1.$ Then
$$(H_{3}f_{1})(t)=1+\frac{1}{2}\cdot\sqrt[7]{4\left(t-\frac{1}{2}\right)}\cdot\int^{\frac{1}{2}}_{-\frac{1}{2}}u^{\frac{1}{7}}du=1=f_{1}(t), \,\ t\in[0,1].$$
b) We define the function $f_{2}$:

$$f_{2}(t)=\frac{1}{2}\left(\sqrt{\frac{57}{17}}+\sqrt{\frac{33}{119}}\cdot\sqrt[7]{2\left(t-\frac{1}{2}\right)}\right), \,\ t\in[0,1].$$
Then $f_{2}\in C[0,1]$ and the function $f_{2}(t)$ is strictly
positive. Put $$a=\frac{1}{2}\sqrt{\frac{57}{17}} \,\ , \,\
b=\frac{1}{2}\sqrt{\frac{33}{119}} \,\ .$$ We have
$$(H_{3}f_{2})(t)=h_{1}(t)+h_{2}(t)+h_{3}(t)+h_{4}(t)+\gamma \,\ ,$$
where
$$h_{1}(t)=\frac{a^{3}}{2}
\varphi(t)\cdot\int^{1}_{0}\sqrt[7]{u-\frac{1}{2}}du,$$

$$h_{2}(t)=\frac{3a^{2}b}{2}\cdot\sqrt[7]{2}
\varphi(t)\cdot\int^{1}_{0}\sqrt[7]{\left(u-\frac{1}{2}\right)^{2}}du,$$

$$h_{3}(t)=\frac{3ab^{2}}{2}\cdot\sqrt[7]{4}
\varphi(t)\cdot\int^{1}_{0}\sqrt[7]{\left(u-\frac{1}{2}\right)^{3}}du,$$

$$h_{4}(t)=\frac{b^{3}}{2}\cdot\sqrt[7]{8}
\varphi(t)\cdot\int^{1}_{0}\sqrt[7]{\left(u-\frac{1}{2}\right)^{4}}du,$$

$$\gamma=\int^{1}_{0}f^{3}_{2}(u)du.$$
Here $\varphi(t)=\sqrt[7]{4\left(t-\frac{1}{2}\right)} \,\ , \,\
t\in[0,1]. $

It is clear that
$$h_{1}(t)=h_{3}(t)\equiv0.$$

For the functions $h_{2}(t)$ and $h_{4}(t)$ we obtain, that

$$h_{2}(t)=\frac{3a^{2}b\sqrt[7]{2}}{2}\cdot\varphi(t)\int^{\frac{1}{2}}_{-\frac{1}{2}}
u^{\frac{2}{7}}du=\frac{7a^{2}b}{6\sqrt[7]{2}}\cdot\varphi(t),$$

$$h_{4}(t)=\frac{b^{3}\sqrt[7]{8}}{2}\cdot\varphi(t)\int^{\frac{1}{2}}_{-\frac{1}{2}}
u^{\frac{4}{7}}du=\frac{7b^{3}}{22\sqrt[7]{2}}\cdot\varphi(t).$$

Observe that
$$\gamma=a^{3}+3ab^{2}\sqrt[7]{4}\cdot\int^{\frac{1}{2}}_{-\frac{1}{2}}u^{\frac{2}{7}}du
=a^{3}+\frac{7ab^{2}}{3}=a.$$
Consequently, we have
$$H_{3}f_{2}=h_{2}+h_{4}+a=a+\frac{7b}{2\sqrt[7]{2}}\left(\frac{a^2}{3}+
\frac{b^{2}}{11}\right)\varphi(t)=a+\frac{b}{\sqrt[7]{2}}\varphi(t)=f_{2}(t).$$
From Proposition \ref{p4.1.}, Lemma \ref{l2.1.} and Proposition \ref{p1} we get

\begin{thm}\label{t4.2.} The model
$$H(\sigma)=-\frac{1}{\beta}\sum_{<x,y>\atop{x,y\in V}}\ln\left(1+\frac{1}{2}
\sqrt[7]{4\left(\sigma(x)-\frac{1}{2}\right)\left(\sigma(y)-\frac{1}{2}\right)}\right)
\,\ , \,\ \sigma\in\Omega_{V}$$ on the Cayley tree $\Gamma^{3}$
has at least two translation-invariant Gibbs measures.
\end{thm}

\section{Existence of two Gibbs measures for the model (\ref{e1.1}): case $k\geq4$}

Let $k\in\mathbb{N}$ and $k\geq2.$ We consider sequences of
continuous functions $P_{n}(x)\,\ (n\in\mathbb{N})$ and
$Q_{m}(x)\,\ (m\in\mathbb{N}, \,\ m>k)$ defined by

$$P_{n}(x)\equiv P_{n,k}(x)=\left(1+\frac{x^{n-1}}{2}\right)^{k+1}-\left(1-\frac{x^{n-1}}{2}\right)^{k+1}, \,\ x\in\mathbb{R},$$

$$Q_{m}(x)\equiv Q_{m,k}(x)=(k+1)x^{m-k}, \,\ m>k, \,\
x\in\mathbb{R}.$$

\begin{pro}\label{p5.1.} Let $k\geq 2.$ Then
\begin{equation}\label{e5.1}
P_{n}(1)>Q_{n}(1),
\end{equation}
for any $n\in\mathbb{N}, \,\ n>k$.
\end{pro}
\proof Let $k\geq2$ and $n>k.$ We have
$$P_{n}(1)=\mu_{k}=\frac{3^{k+1}-1}{2^{k+1}} \,\ , \,\ Q_{n}=\eta_{k}=k+1.$$
In the case $k=2$ we obtain, that
$$P_{n}(1)=\frac{13}{4}>Q_{n}(1)=3.$$
We now suppose, that the inequality (\ref{e5.1}) holds for $k=m>2.$ Then
we show that the inequality (\ref{e5.1}) also is true for $k=m+1.$

Obviously, that

$$\mu_{m+1}=\frac{3^{{(m+1)}+1}-1}{2^{(m+1)+1}}>\frac{3^{(m+1)+1}-3}{2^{m+1}\cdot2}=
\frac{3^{m+1}-1}{2^{m+1}}\cdot\frac{3}{2}=\mu_{m}\cdot\frac{3}{2}>(m+1)\cdot\frac{3}{2}>m+2=\eta_{m+1},$$
i.e. $\mu_{m+1}>\eta_{m+1}.$ Thus we get

$$P_{n}(1)>Q_{n}(1)$$
for any $k\geq2$ and $n>k.$

\begin{pro}\label{p5.2.} Let $k\geq 2.$ The equation
\begin{equation}\label{e5.2}
\left(1+\frac{x}{2}\right)^{k+1}-\left(1-\frac{x}{2}\right)^{k+1}-(k+1)x=0,\,\ \ x\geq 0
\end{equation}
has a unique solution $x=0.$
\end{pro}
\proof Let $k\geq2.$ Define the continuous function
$\varphi(x):$

$$\varphi(x)=\left(1+\frac{x}{2}\right)^{k+1}-\left(1-\frac{x}{2}\right)^{k+1}-(k+1)x, \,\ x\in[0,\infty).$$

We have
$$\varphi'(x)=(k+1)\left(\frac{1}{2}\left(1+\frac{x}{2}\right)^{k}+\frac{1}{2}\left(1-\frac{x}{2}\right)^{k}-1\right).$$

However,
$$\left(1+\frac{x}{2}\right)^{k}+\left(1-\frac{x}{2}\right)^{k}> 2, \,\ \mbox{for all} \,\
x\in(0,\infty).$$\\
Consequently, we have $\varphi'(x)>0$ for all $x\in(0,\infty),$
i.e. the function $\varphi(x)$ is an increasing on $[0,\infty).$ So,
the zero is a unique solution of the equation (\ref{e5.2}).

\begin{pro}\label{p5.3.} Let $k\geq2.$ Then for each
$n\in\mathbb{N}, \,\ n>k $ the equation
\begin{equation}\label{e5.3}
P_{n}(x)-Q_{n}(x)=0
\end{equation}
has at least one solution $\xi=\xi(k;n)$ in (0,1).
\end{pro}
\proof Let $k\geq2$ and $n>k.$ We have

$$\lim\limits_{x\rightarrow0+} \frac{P_{n}(x)}{Q_{n}(x)}=\frac{1}{k+1}\lim\limits_{x\rightarrow0+}
\frac{\left(1+\frac{x^{n-1}}{2}\right)^{k+1}-\left(1-\frac{x^{n-1}}{2}\right)^{k+1}}{x^{n-k}}=$$

$$=\frac{1}{k+1}\lim\limits_{x\rightarrow0+}\frac{\left(\left(1+\frac{x^{n-1}}{2}
\right)-\left(1-\frac{x^{n-1}}{2}\right)\right)\sum\limits_{j=0}^{k}
\left(1+\frac{x^{n-1}}{2}\right)^{k-j}\left(1-\frac{x^{n-1}}{2}\right)^{j}}{x^{n-k}}=$$
$$=\frac{1}{k+1}\lim\limits_{x\rightarrow0+}x^{k-1}\cdot\sum\limits_{j=o}^{k}
\left(1+\frac{x^{n-1}}{2}\right)^{k-j}\left(1-\frac{x^{n-1}}{2}\right)^{j}=0.$$
Since the functions $P_{n}(x)$ and $Q_{n}(x)$ are continuous, the
exists a number $\delta>0$ such that
$$P_{n}(x)<Q_{n}(x) \,\ \mbox{for all} \,\ x\in(0,\delta).$$
However $P_{n}(0)=Q_{n}(0)=0$ and by Proposition \ref{p5.1.} we have
$P_{n}(1)>Q_{n}(1).$ Consequently, there exists a number
$\xi=\xi(k;n)\in(0,1)$ such that
$P_{n}(\xi(k;n))=Q_{n}(\xi(k;n))=0.$

Let $k\geq2$ be a fixed number and suppose that
$\left\{\xi(k;n)\right\}_{n>k}\subset(0,1)$ -- some set of
solutions of the following system of equations:
$$P_{n}(x)-Q_{n}(x)=0, \,\ n\in\mathbb{N}, \,\ n>k.$$
 We have $0<\xi(k;n)<1$ for all $n\in\mathbb{N}, \,\ n>k.$
 Consequently $0<\xi(k;n)^{n-1}<1$ for all $n>k.$ Then there
 exists a upper limit of the sequence $\xi(k;n)^{n-1}, \,\ n>k,$
 i.e. there exists a subsequence $\alpha_{p}=\xi(k;n_{p})^{n_{p}-1}, \,\
 p\in\mathbb{N}$ of the sequence $\xi(k;n)^{n-1}, \,\ n>k$ such
 that
 $$\alpha=\lim\limits_{n\rightarrow\infty}\sup\xi(k;n)^{n-1}=
 \lim\limits_{p\rightarrow\infty}\xi(k;n_{p})^{n_{p}-1}=\lim\limits_{p\rightarrow\infty}\alpha_{p}.$$
 Obviously, that $0\leq\alpha\leq1.$ Define the sequence $\beta_{p}, \,\
 p\in\mathbb{N}$ by

 $$\beta_{p}=\xi(k;n_{p}), \,\ p\in\mathbb{N}.$$
 Then
 $$\alpha_{p}=\beta^{n_{p}-1}_{p}, \,\ p\in\mathbb{N}.$$

\begin{lemma}\label{l5.4.}
 $\alpha=\lim\limits_{p\rightarrow\infty}\alpha_{p}=0.$
 \end{lemma}
\proof a)Assume $\alpha=1.$ Put
 $$\beta=\lim\limits_{p\rightarrow\infty} \sup\xi(k;n_{p})=\lim\limits_{p\rightarrow\infty} \sup\beta_{p}.$$
 Then, there exists a subsequence
 $\{\beta_{p_{q}}\}_{q\in\mathbb{N}}\subset\{\beta_{p}\}_{p\in\mathbb{N}}$
 such that $$\lim\limits_{q\rightarrow\infty}\beta_{p_{q}}=\beta.$$
 We have
 $0\leq\beta\leq1$. If $0\leq\beta< 1$, there exists
 $q_{0}\in\mathbb{N}$ such that $\beta_{p_{q}}<\frac{1+\beta}{2}$
 for all $q>q_{0}.$ From that
 $$0\leq\alpha_{p_{q}}\leq\left(\frac{1+\beta}{2}\right)^{n_{p_{q}}-1}, \,\ q\in\mathbb{N}, \,\ q>q_{0}.$$
 Therefore
 $\alpha=\lim\limits_{q\rightarrow\infty}\alpha_{p_{q}}=0.$
 The last equality is a contradiction to the assumption $\alpha=1$. 
 However, we obtain that $\beta=1.$ Then from the equality
 \begin{equation}\label{e5.4}
  P_{n_{p_{q}}}\left(\xi(k;n_{p_{q}})\right)=
  Q_{n_{p_{q}}}\left(\xi(k;n_{p_{q}})\right), \,\
  q\in\mathbb{N}
  \end{equation}
   as $q\rightarrow\infty$ we observe that
 $$\left(1+\frac{1}{2}\right)^{k+1}-\left(1-\frac{1}{2}\right)^{k+1}=k+1,$$
 i.e. $$P_{m}(1)=Q_{m}(1), \,\ m>k.$$ The last equality is a
 contradiction to the assertion of Proposition \ref{p5.1.}. Thus, we have
 proved that $\alpha\neq1.$

 b) Assume that $0<\alpha<1.$ In the case $0\leq\beta<1$ we get
 $\alpha=0.$ So $\beta=1.$ Then from (\ref{e5.4}) as $q\rightarrow\infty$
 we get
 $$\left(1+\frac{\alpha}{2}\right)^{k+1}-\left(1-\frac{\alpha}{2}\right)^{k+1}=(k+1)\alpha.$$
 The last equality is contradict to the assertion of Proposition \ref{p5.2.}.
 Thus, we have proved that $\alpha \not \in (0,1).$
 Consequently, $\alpha=0.$\endproof

\begin{cor}\label{c5.5.}
 $\lim\limits_{p\rightarrow\infty}\beta_{p}=1.$
 \end{cor}
\proof From the equality (\ref{e5.4}) we get

 $$\beta_{p}=\xi(k;n_{p})=\sqrt[k-1]{\frac{k+1}{\sum\limits_{j=0}^{k}
 \left(1+\frac{\alpha_{p}}{2}\right)^{k-j}\left(1-\frac{\alpha_{p}}{2}\right)^{j}}} \,\ , \,\ p\in\mathbb{N}.$$

Hence by Lemma \ref{l5.4.} it follows that
$$\lim\limits_{p\rightarrow\infty}\beta_{p}=1.$$
\endproof
Define the sequence
$C_{n}, \,\ n>k\geq 2:$
\begin{equation}\label{e5.5}
C_{n}=C_{n}(k)=\frac{\xi(k;n)^{3n-k-2}}{\frac{1}{2+k}\cdot\left[\left(1+
\frac{\xi(k;n)^{n-1}}{2}\right)^{k+2}-
\left(1-\frac{\xi(k;n)^{n-1}}{2}\right)^{k+2}\right]-\xi(k;n)^{n-k}},
\end{equation}
where $\xi(k;n)\in(0,1)$  is an arbitrary solution to the
equation (\ref{e5.3}).

Put $$\gamma_{p}=\gamma_{p}(k)=C_{n_{p}}(k), \,\ p\in\mathbb{N}.$$

\begin{lemma}\label{l5.6.} For every $k\in\mathbb{N}, \,\
k\geq2$ the following equality holds
$$\lim\limits_{p\rightarrow\infty}\gamma_{p}(k)=\frac{12}{k}.$$
\end{lemma}
\proof We have
$$\gamma_{p}=\frac{\alpha^{3}_{p}\cdot\beta^{1-k}_{p}}
{\frac{1}{k+2}\cdot\left(\left(1+\frac{\alpha_{p}}{2}\right)^{k+2}-
\left(1-\frac{\alpha_{p}}{2}\right)^{k+2}\right)-\xi(k;n_{p})^{n_{p}-k}}=$$
$$=\frac{\alpha^{3}_{p}\cdot\beta^{1-k}_{p}}
{\frac{1}{k+2}\cdot\left(\left(1+\frac{\alpha_{p}}{2}\right)^{k+2}-\left(1-\frac{\alpha_{p}}{2}\right)^{k+2}\right)-
\frac{1}{k+1}\cdot\left(\left(1+\frac{\alpha_{p}}{2}\right)^{k+1}-\left(1-\frac{\alpha_{p}}{2}\right)^{k+1}\right)}
\,\ .$$ However
$$\left(1+\frac{\alpha_{p}}{2}\right)^{k+2}-\left(1-\frac{\alpha_{p}}{2}\right)^{k+2}=
\sum\limits_{j=0}^{k+2}C^{j}_{k+2}\cdot\left(\frac{\alpha_{p}}{2}\right)^{j}-
\sum\limits_{j=0}^{k+2}C^{j}_{k+2}\cdot\left(-\frac{\alpha_{p}}{2}\right)^{j}=$$
$$=2C^{1}_{k+2}\cdot\frac{\alpha_{p}}{2}+2C^{3}_{k+2}\cdot\frac{\alpha^{3}_{p}}{2^{3}}+
... +2C^{m_{1}}_{k+2}\cdot\frac{\alpha^{m_{1}}_{p}}{2^{m_{1}}} \,\
,$$ where
$$m_{1}\equiv m_{1}(k)=\left\{%
\begin{array}{ccc}
  k+2, & \mbox{if} \,\ k \,\  \mbox{is odd} \\
  k+1, & \mbox{if} \,\ k \,\  \mbox{is even}. \\
\end{array}\right.
$$ Analogously we have

$$\left(1+\frac{\alpha_{p}}{2}\right)^{k+1}-\left(1-\frac{\alpha_{p}}{2}\right)^{k+1}
=2C^{1}_{k+1}\cdot\frac{\alpha_{p}}{2}+2C^{3}_{k+1}\cdot\frac{\alpha^{3}_{p}}{2^{3}}+
... +2C^{m_{2}}_{k+1}\cdot\frac{\alpha^{m_{2}}_{p}}{2^{m_{2}}},$$
where
$$m_{2}\equiv m_{2}(k)=\left\{
\begin{array}{ccc}
  k+1, & \mbox{if} \,\ k \,\  \mbox{is even} \\
  k, & \mbox{if} \,\ k \,\  \mbox{is odd}, \\
\end{array}\right.
$$ i.e. $m_{2}=2m_{0}-1, \,\ m_{0}\in\mathbb{N}.$

Therefore
$$
{\frac{1}{k+2}\cdot\left(\left(1+\frac{\alpha_{p}}{2}\right)^{k+2}-\left(1-\frac{\alpha_{p}}{2}\right)^{k+2}\right)-
\frac{1}{k+1}\cdot\left(\left(1+\frac{\alpha_{p}}{2}\right)^{k+1}-\left(1-\frac{\alpha_{p}}{2}\right)^{k+1}\right)}=$$
$$=\sum\limits_{j=2}^{m_{0}}a_{j}\alpha_{p}^{2j-1}+a_{m_{0}+1}\alpha_{p}^{2m_{0}+1}=
\alpha_{p}^{3}(a_{2}+a_{3}\alpha_{p}^{2}+a_{4}\alpha_{p}^{4}+\ldots+a_{m_{0}+1}\alpha_{p}^{2(m_{0}-1)}),$$
where
$$a_{j}=\frac{2}{2^{2j-1}}\cdot\left(\frac{C^{2j-1}_{k+2}}{k+2}-\frac{C^{2j-1}_{k+1}}{k+1}\right), \,\ j=2,3,... \,\
,$$
$$a_{m_{0}+1}=\left\{
\begin{array}{ccc}
  0 & \mbox{if} \,\ m_{1}=m_{2}, \\
  \frac{1}{2^{2m_{0}}}\cdot\frac{C^{2m_{0}+1}_{k+2}}{k+2} & \mbox{if} \,\ m_{2}<m_{1}. \\
\end{array}\right.
$$ Obviously that $$a_{2}=\frac{k}{12}.$$

Thus we get

$$\gamma_{p}=\frac{\beta_{p}^{1-k}}{\frac{k}{12}+a_{3}\alpha_{p}^{2}+a_{4}\alpha_{p}^{4}+...+a_{m_{0}+1}\alpha_{p}^{2(m_{0}-1)}},
  \,\ p\in\mathbb{N}.$$

Hence by Corollary \ref{c5.5.} it follows that

$$\lim\limits_{p\rightarrow\infty} \gamma_{p}=\frac{12}{k}\,\
.$$

\begin{cor}\label{c5.7.} If $k\geq4$ then
$0<\lim\limits_{p\rightarrow\infty}\gamma_{p}\leq3.$
\end{cor}

For each $k\geq4$ we define the set
$\mathbb{N}_{0}(k):$
$$\mathbb{N}_{0}(k)=\{p\in\mathbb{N}: |\gamma_{p}(k)|<4\}.$$ Note that, the set
$\mathbb{N}_{0}(k)$ is a countable subset in the set of all
natural numbers. For each $p\in\mathbb{N}_{0}(k), \,\ (k\geq4)$ we
define the continuous function $ K_{p}(t,u;k) $ on $[0,1]^{2}$ by
$$
K_{p}(t,u;k)=1+\gamma_{p}(k)\left(t-\frac{1}{2}\right)\left(u-\frac{1}{2}\right),
\,\ t,u\in[0,1].$$ By the inequality $|\gamma_{p}(k)|<4$ it
follows that, the
function $K_{p}(t,u;k)$ is strictly positive.

\begin{thm}\label{t5.8.} Let $k\geq4.$ For each
$p\in\mathbb{N}_{0}(k)$ the Hammerstein's equation
\begin{equation}\label{e5.6}
\int^{1}_{0}K_{p}(t,u;k)f^{k}(u)du=f(t)
\end{equation}
in the $C[0,1]$ has at least two positive solutions.
\end{thm}
\proof Obviously, that the function $f_{0}(t)\equiv1$ is a solution of
the equation (\ref{e5.6}). Define the strictly positive continuous
function $f_{1}(t)$ on $[0,1]$ by
$$
f_{1}(t)=\xi(k;n_{p})+\xi(k;n_{p})^{n_{p}}\left(t-\frac{1}{2}\right),\,\
t\in[0,1].$$

We shall prove that the function $f_{1}(t)$ also is a solution of
the Hammerstein's equation $(\ref{e5.6}):$

$$\int_{0}^{1}K_{p}(t,u;k)f_{1}^{k}(u)du=
\int^{1}_{0}\left(1+\gamma_{p}(k)\left(t-\frac{1}{2}\right)\left(u-\frac{1}{2}\right)\right)\times$$
$$\times\left(\xi(k;n_{p})+\xi(k;n_{p})^{n_{p}}\left(u-\frac{1}{2}\right)\right)^{k}du=
\int^{1/2}_{-1/2}\left(1+\gamma_{p}(k)\left(t-\frac{1}{2}\right)u\right)\left(\beta_{p}+\beta^{n_{p}}_{p}u\right)^{k}du=$$
$$=\int^{1/2}_{-1/2}\left(\beta_{p}+\beta^{n_{p}}_{p}u\right)^{k}du+\gamma_{p}(k)
\left(t-\frac{1}{2}\right)\int^{1/2}_{-1/2}u\left(\beta_{p}+\beta^{n_{p}}_{p}u\right)^{k}du=$$
$$=\frac{\beta^{k}_{p}}{\beta^{n_{p}-1}_{p}}\int^{1/2}_{-1/2}\left(1+\beta^{n_{p}-1}_{p}u\right)^{k}
d\left(1+\beta^{n_{p}-1}_{p}u\right)+\gamma_{p}(k)\left(t-\frac{1}{2}\right)\times$$
$$\times\frac{\beta^{k}_{p}}{\beta^{n_{p}-1}_{p}}\int^{1/2}_{-1/2}
u\left(1+\beta^{n_{p}-1}_{p}u\right)d\left(1+\beta^{n_{p}-1}_{p}u\right)=
\frac{\beta^{k}_{p}}{\alpha_{p}}\cdot\frac{1}{k+1}\left(1+\alpha_{p}u\right)^{k+1}\left|^{1/2}_{-1/2}+\right.$$
$$+\frac{\gamma_{p}(k)\beta^{k}_{p}}{\alpha^{2}_{p}}\left(t-\frac{1}{2}\right)\int^{1/2}_{-1/2}
\left(\left(1+\alpha_{p}u\right)^{k+1}-\left(1+\alpha_{p}u\right)^{k}\right)d(1+\alpha_{p}u)=$$
$$=\frac{\beta^{k}_{p}}{\alpha_{p}}\cdot\frac{1}{k+1}\left(\left(1+\frac{\alpha_{p}}{2}\right)^{k+1}-
\left(1-\frac{\alpha_{p}}{2}\right)^{k+1}\right)+\frac{\gamma_{p}(k)\beta^{k}_{p}}{\alpha^{2}_{p}}\times
\left(t-\frac{1}{2}\right)\cdot\left(\frac{1}{k+2}\left(1+\alpha_{p}u\right)^{k+2}\left|^{1/2}_{-1/2}-\right.\right.$$
$$\left.-\frac{1}{k+1}\left(1+\alpha_{p}u\right)^{k+1}|^{1/2}_{-1/2}\right)=\frac{\beta^{k}_{p}}{\alpha_{p}}\cdot
\frac{1}{k+1}\cdot(k+1)\beta^{n_{p}-k}_{p}+\frac{\gamma_{p}(k)\beta^{k}_{p}}{\alpha^{2}_{p}}
\cdot\left(t-\frac{1}{2}\right)\times$$
$$\times\left[\frac{1}{k+2}\left(\left(1+\frac{\alpha_{p}}{2}\right)^{k+2}-
\left(1-\frac{\alpha_{p}}{2}\right)^{k+2}\right)-\frac{1}{k+1}\left(\left(1+\frac{\alpha_{p}}{2}\right)^{k+1}-
\left(1-\frac{\alpha_{p}}{2}\right)^{k+1}\right)\right]=$$
$$=\frac{\beta^{n_{p}}_{p}}{\alpha_{p}}+\frac{\gamma_{p}(k)\beta^{k}_{p}}{\alpha^{2}_{p}}\cdot\left(t-\frac{1}{2}\right)\cdot
\left[\frac{1}{k+2}\left(\left(1+\frac{\alpha_{p}}{2}\right)^{k+2}-
\left(1-\frac{\alpha_{p}}{2}\right)^{k+2}\right)-\frac{1}{k+1}(k+1)\beta^{n_{p}-k}_{p}\right]=$$
$$=\beta_{p}+\frac{\gamma_{p}(k)\beta^{k}_{p}}{\alpha^{2}_{p}}\cdot\left(t-\frac{1}{2}\right)
\left(\frac{1}{k+2}\left(\left(1+\frac{\alpha_{p}}{2}\right)^{k+2}-
\left(1-\frac{\alpha_{p}}{2}\right)^{k+2}\right)-\beta^{n_{p}-k}_{p}\right)=$$
$$=\beta_{p}+\frac{\beta^{k}_{p}}{\alpha^{2}_{p}}\cdot
\alpha^{3}_{p}\beta^{1-k}_{p}\left(t-\frac{1}{2}\right)=\beta_{p}+\beta^{n_{p}}_{p}\cdot
\left(t-\frac{1}{2}\right)=\xi(k;n_{p})+\xi(k;n_{p})^{n_{p}}\left(t-\frac{1}{2}\right)=f_{1}(t).$$
\endproof

From Theorem \ref{t5.8.}, Lemma \ref{l2.1.} and Proposition \ref{p1} we get the following theorem.

\begin{thm}\label{t5.9.} Let $k\geq4$
and $p\in\mathbb{N}_{0}(k).$ The model
$$H(\sigma)=-{1\over\beta}\sum\limits_{<x,y>\atop{x,y}\in
V}\ln K_{p}(\sigma(x),\sigma(y);k), \,\,\,\
\sigma\in\Omega_V$$
on the Cayley tree $\Gamma^{k}$ has at
least two translations-invariant Gibbs measures.
\end{thm}
\section*{ Acknowledgements}

 U. Rozikov thanks Institut des Hautes \'Etudes Scientifiques (IHES), Bures-sur-Yvette, France for support of his visit to IHES and IMU/CDC-program for a (travel) support.

\end{document}